\documentclass[final,5p,times,authoryear]{elsarticle} 

\usepackage{graphicx}
\usepackage{amsmath}

\usepackage{accents}   

\usepackage{array}

\usepackage{url}

\usepackage{nomencl}
\usepackage{mathtools}

\usepackage{lineno} 
\usepackage{hyperref}
\usepackage{subcaption}
\usepackage{booktabs}
\usepackage{bm}
\usepackage{algorithmicx}
\usepackage{amssymb}
\usepackage{diagbox}

\usepackage{amsthm}

\theoremstyle{plain}

\modulolinenumbers[5]

\journal{arXiv}

\biboptions{numbers,sort&compress}

\begin{document}

\begin{frontmatter}

\title{Extreme events in time series aggregation: A case study for optimal residential energy supply systems}

\author[ere]{Holger Teichgraeber \tnoteref{t1}}
\ead{hteich@stanford.edu}
\author[ere]{Constantin P. Lindenmeyer \tnoteref{t1}}
\ead{constantin.lindenmeyer@rwth-aachen.de}
\author[ltt]{Nils Baumg\"artner}
\author[fzj]{Leander Kotzur}
\author[fzj,cfc]{Detlef~Stolten}
\author[fzj]{Martin Robinius}
\author[ltt,jul]{Andr\'{e} Bardow}
\author[ere]{Adam R. Brandt \corref{cor}}
\ead{abrandt@stanford.edu}

\cortext[cor]{Corresponding author. Tel: 650-724-8251}

\tnotetext[t1]{HT and CPL contributed equally to this work.}

\address[ere]{Department of Energy Resources Engineering, Stanford University,
  Green Earth Sciences Building 065, 367 Panama St., Stanford, California, USA}
  
\address[ltt]{Institute of Technical Thermodynamics, RWTH Aachen University, Schinkelstrasse 8, 52062 Aachen, Germany}

\address[fzj]{Institute of Energy and Climate Research - Techno-economic System Analysis, Forschungszentrum J\"ulich GmbH, Wilhelm-Johnen-Str. 52428 J\"ulich, Germany}

\address[jul]{Institute of Energy and Climate Research - Energy Systems Engineering (IEK-10), Forschungszentrum J\"ulich GmbH, 52428 Jü\"ulich, Germany}

\address[cfc]{Chair for Fuel Cells, RWTH Aachen University, c/o Institute of Electrochemical Process Engineering (IEK-3), Forschungszentrum J\"ulich GmbH, Wilhelm-Johnen-Str., 52428 J\"ulich, Germany}
  
\begin{abstract}

To account for volatile renewable energy supply, energy systems optimization problems require high temporal resolution. Many models use time-series clustering to find representative periods to reduce the amount of time-series input data and make the optimization problem computationally tractable. However, clustering methods remove peaks and other extreme events, which are important to achieve robust system designs. 
We present a general decision framework to include extreme events in a set of representative periods. 
We introduce a method to find extreme periods based on the slack variables of the optimization problem itself. 
Our method is evaluated and benchmarked with other extreme period inclusion methods from the literature for a design and operations optimization problem: a residential energy supply system.  
Our method ensures feasibility over the full input data of the residential energy supply system although the design optimization is performed on the reduced data set.

We show that using extreme periods as part of representative periods improves the accuracy of the optimization results by 3\% to more than 75\% depending on system constraints compared to results with clustering only, and thus reduces system cost and enhances system reliability.

\end{abstract}
  

\end{frontmatter}


\makenomenclature 

\section{Introduction}

In order to comply with the Paris Agreement, greatly increasing penetration rates of renewable energies will be needed. For example, the 2030 Energy Strategy of the European Union aims at a share of renewable energies in the energy sector of more than 32\% by 2030. Ideally, these future energy systems will be thoroughly optimized to ensure least cost and high reliability. For example, capacity expansion planning (CEP) optimization models optimize future grid generation mixes and grid topologies to reduce emissions at least cost \citep{Merrick2016a}. On a more detailed level, residential clean energy systems optimization can be performed to reduce the cost of meeting stringent clean energy targets while minimizing call on the grid \citep{Kotzur2018}.

A challenging feature of such models is that long-term decisions with effective scales of years or decades (e.g., capital investment in wind farm) have to be merged with operations optimization based on temporal weather and demand data (e.g., hourly wind output at a given site) \citep{Despres2015, Nahmmacher2016}. Thus, the design optimization problem for such kind of volatile and decentralized energy systems are computationally expensive \citep{Goderbauer2019TheNP-hard}, because of the size of the necessary input data when the design must co-consider high frequency temporal data (wind speed, solar availability etc.) \citep{Mancarella2014a, Bahl2017a, Kotzur2018}.   Even more challenging are systems with large geographic scale: to properly model such system, high resolution in space and time is necessary, which lead to even larger input data sets and more computational challenge to solve \citep{Nahmmacher2016, Kotzur2018}. 

In order to build optimization models that are computationally tractable, the complexity of the problem needs to be reduced. The complexity of such a problem is reduced in three general ways: reducing physical detail \citep{Geidl2007, Banos2011}, spatial aggregation \citep{Mancarella2014a}, and temporal aggregation  \citep{Petruschke2014, Stadler2014}. In practice, any energy systems optimization model uses a combination of some aspect of these three \citep{Collins2017IntegratingReview}. Reducing the level of physical detail is a common approach, such as when the complexities of grid-level AC power flow are simplified to a DC optimal power flow problem, which is much easier to solve \citep{Cain2012HistoryFormulations}. In spatial aggregation, adjacent regions are merged into a larger region, reducing model geographic fidelity but also reducing the size of input data and the number decision variables \citep{Pfenninger2014EnergyChallanges, Collins2017IntegratingReview}. Similarly, temporal aggregation replaces the full fidelity time series data with simplified representations, e.g. representative periods \citep{Pfenninger2017DealingVariability}. 

For temporal aggregation, various approaches exist to select representative periods: statistical  and empirical selection processes \citep{Poncelet2015}, graphical aggregation \citep{Ortiga2011}, heuristic selections \citep{Fripp2008, Mavrotas2008, Casisi2009, Nelson2012}, and time-series clustering \citep{Dominguez-Munoz2011, Fazlollahi2014, Green2014, Stadler2014, Adhau2015, Brodrick2015, Merrick2016a, Nahmmacher2016, Bahl2017a, Heuberger2017, Pfenninger2017, Schutz2017, Teichgraeber2017, Bahl2018, Bahl2018a, Gabrielli2018, Kotzur2018, Kotzur2018a, Mallapragada2018, Tejada-Arango2018}. In the last few years, clustering has emerged as a popular method to reduce the number of time-steps used as input data for energy systems optimization problems \citep{Teichgraeber}. Clustering has been applied to input data for multiple optimization problems: large generation capacity expansion problems \citep{Green2014, Merrick2016a, Nahmmacher2016, Heuberger2017, Pfenninger2017, Almaimouni2018, Mallapragada2018}, industrial or residential energy supply system problems \citep{Bahl2017a, Schutz2017, Kotzur2018, Kotzur2018a}, and individual technology optimization problems \citep{Brodrick2015, Teichgraeber2017, Brodrick2017a, Brodrick2018}. 

All of the above mentioned studies use different clustering methods. Teichgraeber and Brandt \citep{Teichgraeber} provide a framework that allows for classification and inter-comparison of different clustering methods.
Some studies also compare the performance of different clustering methods and their effect on the objective function of the optimization problem \citep{, Pfenninger2017, Kotzur2018a}. The most commonly used clustering methods are k-means clustering \citep{Fazlollahi2014, Green2014, Adhau2015, Brodrick2015, Bahl2017a, Heuberger2017, Pfenninger2017, Teichgraeber2017, Gabrielli2018, Kotzur2018, Kotzur2018a}, k-medoids clustering \citep{Dominguez-Munoz2011, Stadler2014, Schutz2017, Bahl2018a, Kotzur2018a}, and hierarchical clustering with the medoid as its representation \citep{Merrick2016a, Nahmmacher2016, Pfenninger2017, Kotzur2018a}. 

In some papers that use clustering, the authors state the concern that peak demands and other extreme events are not properly represented because the clustering methods do not select them as representative periods \citep{Gabrielli2018, Pfenninger2017}. Therefore, these peaks and extreme periods have to be added manually to the clustered input data to maintain the original variation of the input data. This addition serves two purposes: to ensure that operating constraints are met by the designed system (e.g., sufficient supply at all hours),\citep{Bahl2017a, Bahl2018, Bahl2018a, Blanford2018, Gabrielli2018, Baumgaertner2019} and to improve the accuracy of the objective function value of the optimization problem \citep{Dominguez-Munoz2011, Pfenninger2017, Kotzur2018}. Without an adequate representation of extreme events, the accuracy of the optimization results may be compromised \citep{Pfenninger2017}. 

The methods used to select extreme events and how to include them differ from paper to paper. The most common method is to add extreme periods as representative days \citep{Dominguez-Munoz2011, Pfenninger2017, Kotzur2018}. 
In order to maintain integrally preserved clustered input data, the extreme periods have to be excluded from the clustering process prior to clustering \citep{Dominguez-Munoz2011}. These extreme periods are individually added to the optimization problem, whereas all other input data are clustered.  
This method is also used by Fazlollahi et al.~\citep{Fazlollahi2014}, Gabrielli \citep{Gabrielli2018}, and Heuberger et al.~\citep{Heuberger2017}. 

Kotzur et al.~\citep{Kotzur2018} introduce methods modifying the clustering process to include extreme events. Their methods represent certain clusters by extreme periods instead of the centroid or medoid. Comparing these methods with the simple addition of an extreme event (as used by Pfenninger \citep{Pfenninger2017} and Dominguez-Munoz et al.~\citep{Dominguez-Munoz2011}), Kotzur et al.~\citep{Kotzur2018} find no significant differences: Neither the objective function value of the optimization problem nor the individual design variables are significantly different, with all reaching margins of error of below 3\% for the total system cost of a residential energy supply system \citep{Kotzur2018}. 

Bahl et al.~and Baumg\"artner et al.~\citep{Bahl2017a, Bahl2018, Bahl2018a, Baumgaertner2019} use so-called ``feasibility steps'': Extreme events are added to the set of constraints with weight of zero in the objective function of the optimization problem. Feasibility steps ensure that the resulting design from the clustered optimization can supply their system over the entire time of input data. For example, a day with very low wind and solar output will be added to the constraint set to ensure that the resulting design can meet demand on that day. Their approach is iterative: if the operations optimization with the clustered design is infeasible, Bahl et al.~\citep{Bahl2017a} add extreme events as constraints with zero weight in the objective function of the optimization problem and rerun their process until the operations optimization becomes feasible for every day of the full input data.

The used methods in these papers are often not defined or introduced in detail, which leaves room for interpretation. In this paper, we thus first present a general framework to categorize the modeling decisions to be made when using extreme periods. This allows for intercomparability of past and future research.
Second, we introduce a selection method to include extreme periods energy systems optimization problems based on slack variables from the optimization problem itself. 
Third, we show the importance of adding extreme periods to the set of representative periods for energy systems optimization when using clustering.
Finally, we compare existing extreme period selection methods to our new slack-variable-based method on a case study that optimizes the design and operations of a residential energy supply system. 

The structure of this paper is as follows: In Section \ref{sec:Methods}, we describe our methods, including our framework to identify extreme periods. In Section \ref{sec:Data}, we introduce the experimental setup and residential energy supply system optimization problem on which we test the different extreme period selection methods. In Section \ref{sec:Results}, we present and discuss the results, followed by the conclusion in Section \ref{sec:Conclusion}.  

\begin{figure}[htb]
\centering
{\noindent\includegraphics[width=0.5\textwidth]{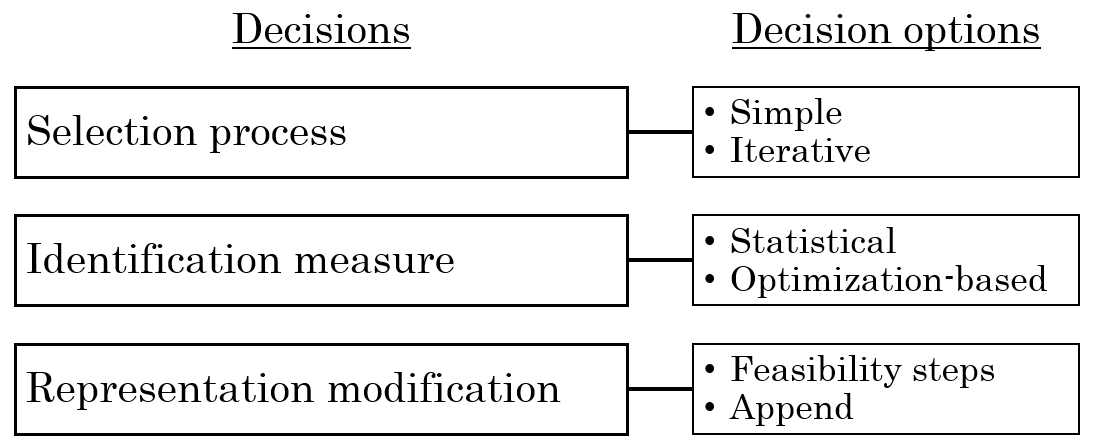}}
\caption[Overall Framework]{Overall framework for the inclusion of extreme
periods when using clustered input data for the optimization of energy systems. Modeling decisions fall into three categories:  the selection process, the identification measure, and the representation modification. In subsections \ref{sec:Selection} to \ref{sec:Representation} we explain each respective decision in detail.}
\label{fig:02_01_Overall_Framework}
\end{figure}

\section{Methods} 
\label{sec:Methods}

Clustering methods to find representative days for the optimization of energy systems have been used in many applications. Teichgraeber and Brandt~\citep{Teichgraeber} introduce a framework for the decisions which have to be made when using clustering. However, the framework does not describe the inclusion of extreme periods as representative periods. 
In Figure~\ref{fig:02_01_Overall_Framework}, we introduce a framework for adding extreme periods to clustered input data. 
Overall, addition of extreme periods requires three decisions that fall into the following categories: (1) the overall selection process for extreme periods, (2) the identification measure for identifying extreme periods, and (3) the representation modification of the clustered input data, i.e.~how to represent extreme periods in representative periods. 

In the following subsections, we first introduce the normalization and clustering methods used in this work, then provide the optimization problem variations used in the extreme period framework, and then outline each step of the extreme period framework. We then present the specific extreme period inclusion methods that we use in this paper.

\subsection{Clustering: Normalization, Algorithm, Representation}
\label{sec:Clust_Algo}
In the description of the clustering method we use in this work, we follow the framework introduced by Teichgraeber and Brandt \citep{Teichgraeber}.
We normalize all data types (demand, solar availability, etc...) individually using z-normalization, where data is normalized to mean $\mu=0$, standard deviation $\sigma=1$. Individual normalization results in one $\mu_d$ and one $\sigma_d$ per data type $d$ across the entire data series (e.g., one $\mu_d$ per year). 
We use the k-means algorithm, as first introduced by Steinhaus and Macqueen \citep{Macqueen1967}. 
k-means employs a partitional clustering algorithm and minimizes the sum of squared distances (SSD) between the members of each cluster, summed over all clusters. It employs the Euclidean distance (ED) as the distance measure between the original periods and the clustered data centers. We use the cluster centroid both during the algorithm computation and to represent the cluster.
Partitional clustering algorithms are greedy and thus, the k-means algorithm only converges to local solutions. We repeat the process with 10,000 initializations and pick the clustering with minimal SSD across all initializations (this allows for reproducibility of results).  Note that choosing the clustering with the lowest SSD between the cluster centers and the original input data does not necessarily result in the best estimate of the objective function of the optimization problem with the full input data \citep{Teichgraeber}. However, this choice does not affect this study, as we compare the k-means and k-means plus extreme period formulations in the domain of the objective function. To perform the clustering, we  use the open-source \href{https://github.com/holgerteichgraeber/TimeSeriesClustering.jl}{TimeSeriesClustering} package \citep{Teichgraeber2019TimeSeriesClustering:Julia}.

\nomenclature{SSD}{Sum of squared distances} 
\nomenclature{$\mu$}{mean}
\nomenclature{$\sigma$}{standard deviation}
\nomenclature{$d$}{input data set}

\subsection{Optimization problem variations}
\label{sec:Opt_Var}

We are interested in approximating the optimization problem with full input data by an optimization problem with a small number of representative days. In the process of selecting and evaluation clustering and extreme periods, we use four different formulations of the same optimization problem based on the different input data sizes. Table \ref{tab:02_01-Overview_Optimization} shows an overview of the four different optimization problem formulations with different input data (clustered or full input data sets), which decision variables are optimized and which are fixed in each problem formulation, and the notation used throughout this paper.

\begin{table*}[h]
\centering
\caption{Overview of different optimization problem formulations reference case, initial optimization, full operations optimization and daily operations optimization, their corresponding number of input periods, decision variables and output notation. \(N_{full}\) stands for the full input data set, \(k\) stands for the number of periods obtained from clustering, \(X\) stands for the number of selected extreme periods, \(DV\) stands for design decision variable, \(OV\) stands for operational decision variable, \(f(x)\) stands for the objective function value of the optimization problem and the asterisk represents an optimal solution.}
\label{tab:02_01-Overview_Optimization}
\begin{tabular}{llllll}
\toprule
\begin{tabular}[t]{@{}l@{}} \textbf{Optimization} \\ \textbf{problem} \\ \textbf{formulation} \end{tabular}& \textbf{Symbol} & 	\begin{tabular}[t]{@{}l@{}} \textbf{Number} \\ \textbf{of} \\ \textbf{periods} \end{tabular}	&	\begin{tabular}[t]{@{}l@{}} \textbf{Decision} \\ \textbf{variables} \end{tabular}	&	\begin{tabular}[t]{@{}l@{}} \textbf{Fixed} \\ \textbf{decision} \\ \textbf{variables} \end{tabular}	&	\begin{tabular}[t]{@{}l@{}} \textbf{Output} \\ \textbf{notation} \end{tabular}	\\

\midrule
Reference case	&	\(O_{ref}\)	&	\(N_{full}\)		&		DV \& OV		&		---		&		\begin{tabular}[t]{@{}l@{}} \(DV_{ref}\) \& \\ \(f(x^{*})_{ref}\) \end{tabular} 		\\

\begin{tabular}[t]{@{}l@{}} Initial \\ optimization \end{tabular}	&	\(O_{init}\)	&	k + X		&		DV \& OV		&		---		&		\begin{tabular}[t]{@{}l@{}} \(DV_{repr}\) \& \\ \(f(x^{*})_{clustered}\) \end{tabular} 		\\

\begin{tabular}[t]{@{}l@{}} Full operations \\ optimization \end{tabular}	&	\(O_{op}\)	&	\(N_{full}\)		&		OV		&		\(DV_{repr}\)		&		\begin{tabular}[t]{@{}l@{}}  \(f(x^{*})_{operations}\) \& \\ feasibility or \\ slack variable \end{tabular}		\\

\begin{tabular}[t]{@{}l@{}} Daily operations \\ optimization \end{tabular}	&	\(O_{daily}\)	&	1		&		OV		&		\(DV_{repr}\)		&		\begin{tabular}[t]{@{}l@{}} feasibility \\  \end{tabular}		\\

\bottomrule
\end{tabular}
\end{table*}

First, the reference case (\(O_{ref}\)) is a design and operations optimization where full input data are used to simultaneously optimize the design and operation of the system. If possible, one would always solve $O_{ref}$, as this would ensure that the optimal design also satisfied feasibility. In general, solving $O_{ref}$ is too computationally challenging, and in these cases the results from \(O_{ref}\) are unknown. Because our experimental setup is simple enough, $O_{ref}$ is still computationally tractable and we can compare the results of any clustering approach to those from $O_{ref}$. 

Second, the initial optimization (\(O_{init}\)) is a design and operations optimization employing clustered input data. $O_{init}$ is computationally tractable because it has less operational decision variables than $O_{ref}$ ($k + X \leq N_{full}$). Because the clustered input data are used instead of the full input data, no guarantee can be made of the feasibility or optimality of the resulting design. 

Third, the full operations optimization (\(O_{op}\)) is an operations optimization where we use the design \(DV_{init}\) obtained from $O_{init}$ and optimize the operations over the full input data. This problem is computationally tractable for many applications because the complicating design variables are fixed. If $O_{op}$ optimization is feasible, the design from the clustered optimization is a viable system design. However, if \(DV_{init}\) is infeasible, this design cannot provide demand over the entire input data set. Note that \(O_{op}\) is the same optimization formulation that Bahl et al.~and Baumg\"artner et al.~\citep{Bahl2017a, Bahl2018, Bahl2018a, Baumgaertner2019} use to evaluate their representative periods. 

Fourth, the daily operations optimization ($O_{daily}$) is an operations optimization using \(DV_{init}\), but optimizes each day individually. The information on feasibility in this formulation allows for identification of infeasible days.  $O_{daily}$ is easier to solve than \(O_{op}\) due to the fact that each day is solved individually.

\nomenclature{k}{number of clusters} 
\nomenclature{$DV$}{Design decision variable}
\nomenclature{$OV$}{Operational decision variable}
\nomenclature{${O}_{ref}$}{Reference case}
\nomenclature{${O}_{init}$}{Initial optimization}
\nomenclature{${O}_{op}$}{Operations optimization}
\nomenclature{${O}_{daily}$}{Daily operations optimization}
\nomenclature{${N}_{full}$}{Full input data set (not clustered)}

\subsection{Framework for the inclusion of extreme periods}

The framework for using extreme events for energy systems optimization problems consists of the three decisions to be made---selection process, identification measure, representation modification---which are presented in the subsections below. All three decisions are independent of each other. See Figure~\ref{fig:02_01_Overall_Framework} for an overview of the framework.

Extreme events can be selected as individual extreme values, e.g. extreme hours, or as extreme periods, e.g. extreme days. This choice largely depends on the optimization problem and its characteristics. If storage is part of the optimization problem and can be used in extreme events, it is important to link successive hours, and in that case, it is important to use extreme periods that cover the storage cycle. However, if storage is not a concern, individual extreme values like extreme hours may be sufficient to capture the characteristics relevant to the optimization problem. 
Because most of applications in the literature link successive hours, we generally refer to extreme periods in the following.

\subsubsection{Selection process}
\label{sec:Selection}

The process of selecting extreme periods can be put into two broad categories: simple extreme period selection and iterative extreme period selection. 

In the first category, simple extreme period selection, individual extreme periods can be added before the optimization run. Using this method, a predefined number of extreme periods is added using any of the identification measures described in Section \ref{sec:Identification}. Simple extreme period selection has been used in the majority of the literature that considers extreme periods, and extreme periods usually are identified with statistical identification measures \citep{Dominguez-Munoz2011, Heuberger2017, Pfenninger2017, Gabrielli2018, Kotzur2018}.

In the second category, extreme periods can be selected iteratively, each after a new optimization run, until a convergence criterion is met. Example convergence criteria can be that all periods of the individual daily optimization are feasible, or that the slack variables in the operations optimization (\(O_{op}\)) are zero (see Section \ref{sec:Identification} for details). 
Iterative extreme period selection has previously been used by \citep{Bahl2017a, Bahl2018, Baumgaertner2019} with identification measures that are based on the optimization problem.

\begin{figure*}[htb]
\centering
{\noindent\includegraphics[width=0.7\textwidth]{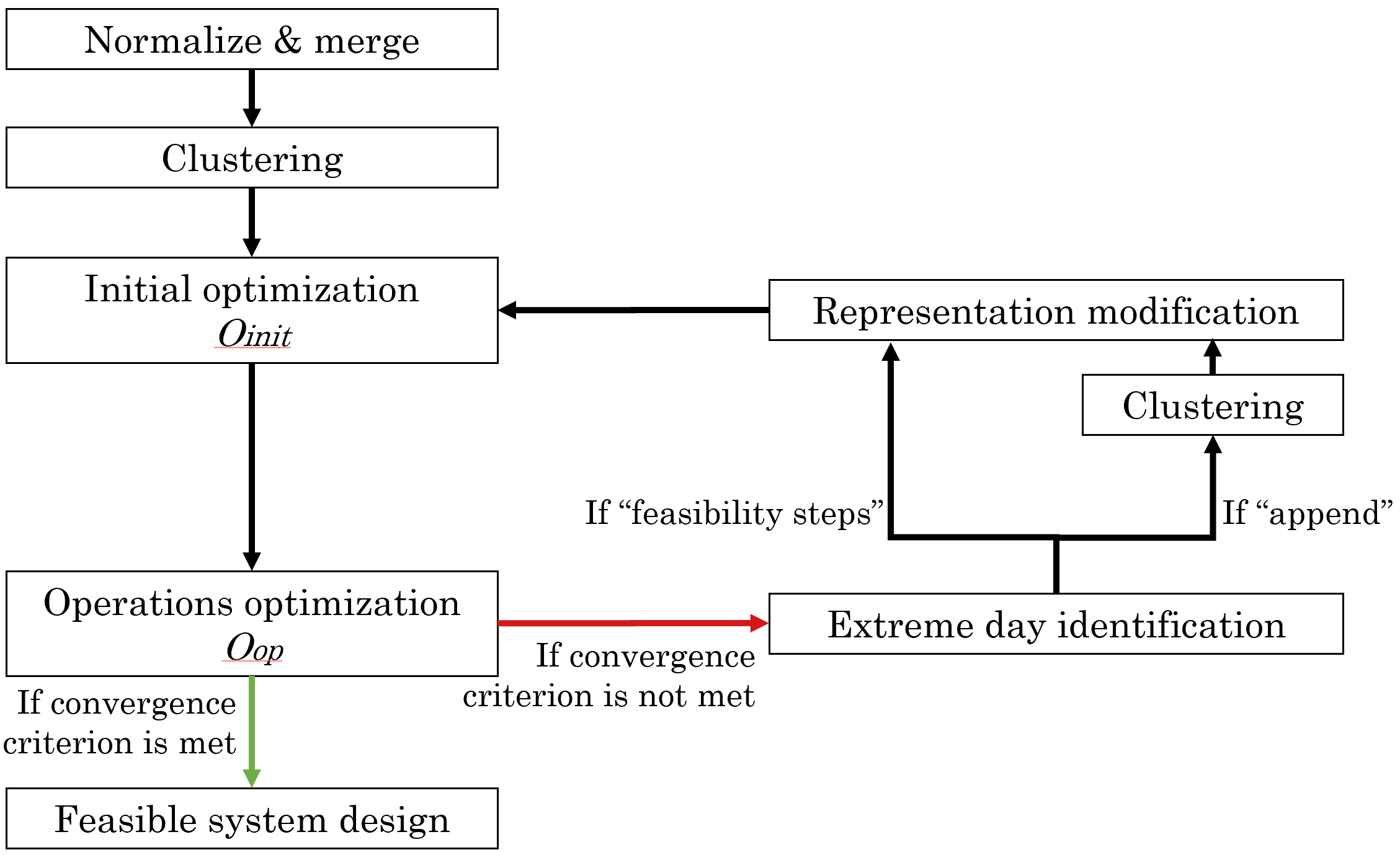}}
\setlength{\abovecaptionskip}{5pt}
\caption[Framework Inclusion]{Overview of the iterative extreme period selection process that uses identification measures from the optimization problem itself.}
\label{fig:02_02_Feasibility_Based}
\end{figure*}

As an example, Figure~\ref{fig:02_02_Feasibility_Based} shows an iterative selection process for convergence criteria based on the optimization problem. The iterative selection process includes the conventional normalization and clustering, but then adds an iterative procedure that includes \(O_{init}\) and \(O_{op}\). If the convergence criterion is not met, a new extreme period is identified based on the identification measure (see Section \ref{sec:Identification}). Then, the extreme periods are added (see Section \ref{sec:Representation}).

In this paper, we call extreme periods that were selected using a simple selection process ``simple'' extreme periods, and we call extreme periods that were selected using an iterative selection process based on the optimization problem ``sufficient'' extreme periods (indicating that they are sufficient for feasibility on the full operations problem (\(O_{op}\))). 

\subsubsection{Identification measure}
\label{sec:Identification}

Extreme periods can be identified using statistical properties of the data or using information from the optimization problem itself. Statistical identification measures work with properties of the data (e.g., take as extreme period the time-step with the highest value of the input data set). Optimization-based identification measures take information of the initial or operations optimization and identify extreme periods based on those results (e.g., day with largest constraint violation). 

Statistical methods identify extreme periods based on extrema (maximum and/or minimum) of the input data sets. Extrema can be absolute extrema or integral extrema. 
An extreme period based on absolute extrema is the period that contains the maximum or minimum absolute value of the overall data within its period. In contrast, an extreme period based on integral extrema is identified as the period that contains the maximum or minimum sum of values compared to all other periods. 
These extrema can be found either for each attribute individually or for all attributes together. If extreme periods based on all attributes together are identified, the data have to be normalized so as to not compare, for example, solar PV insulation in W/m$^2$ to energy prices in \$/kWh. 

One can also identify extreme periods based on information from the optimization problem itself.  The general idea is that extreme periods can be identified by operating the optimal design found in the initial reduced form optimization \(O_{init}\) on the full input data, either as one optimization run \(O_{op}\), or for each day individually \(O_{daily}\) if periods are operationally separable. The identification measures of interest are either the resulting slack variables or feasibility. 

Slack variables that can be used as identification measure are for example virtual generation variables on the demand constraint to provide the unmet energy demand the designed energy system cannot provide. For each energy carrier, one slack variable is added to the demand constraint and can provide virtual power for each time-step if needed, though at a very high cost. A possible extreme period thus can be selected as the period with the highest peak or the highest integral slack variable within one period. For each energy carrier, one slack variable has to be added to the demand constraints and each slack variable needs to be considered for selecting extreme periods. Slack variables have not been previously used as an identification measure in the literature and are newly proposed in this paper.

To use feasibility as an identification measure, one can run an operational optimization \(O_{daily}\) for every period of the original data individually, with the design found in \(O_{init}\). \(O_{daily}\) will be infeasible on some periods. Thus, it is possible to identify which periods are infeasible based on the given system design. Infeasibility identifies these periods as extreme periods. Feasibility-based events were first introduced by Bahl et al.~\citep{Bahl2017a}. Note that feasibility-based identification is only possible if the problem can be separated operationally. In cases where there are constraints that link different periods (e.g. yearly emissions limits or seasonal storage), this approach is not applicable. 

\subsubsection{Representation modification}
\label{sec:Representation}

After identifying the extreme periods, extreme periods have to be added to the clustering representation \citep{Dominguez-Munoz2011, Bahl2017a, Pfenninger2017, Gabrielli2018, Kotzur2018}. We call the methods ``feasibility steps'' and ``append''. The methods differ in the way they change or append the clustering representation and optimization process. 

The ``feasibility steps'' method was introduced by Bahl et al.~\citep{Bahl2017a}. They add an extreme period as a single period cluster to the clustered data. However, this single period cluster is not assigned a weight towards the objective function of the optimization problem, but instead only enters the constraint body. Feasibility steps ensure that the design variables are adequately chosen to handle the extreme period, but do not bias the resulting objective function directly. It is possible to use extreme periods with the same number of time steps as the other periods, or just individual extreme values (e.g. just the absolute maximum heat demand), if the optimization formulation can handle the different formats within the input data arrays \citep{Bahl2017a}. 

The ``append'' method was introduced by Fazlollahi et al.~\citep{Fazlollahi2014}. It finds the chosen extreme period from the input data set, excludes it from the clustering algorithm and afterwards adds this event as a single period cluster to the resulting clusters. For example, in a problem with two input data sets with different extreme periods, the number of periods of the input data would therefore be reduced to \(n_{i} = N_{i} - 2\). After the clustering into \(k\) clusters and the selection of the lowest SSE, these two extreme periods are added to the clustered data as additional representative periods, each with the weight of \(w = 1/N_{i}\). The weights of the other clusters are adjusted accordingly, so that the sum of all \(k+2\) clusters equals 1 again \citep{Kotzur2018}. As long as the centroid instead of the medoid is used as cluster representation, this method preserves the original mean of the data without any additional normalization \citep{Dominguez-Munoz2011}. This method is most likely also used by \citep{Pfenninger2017, Gabrielli2018}, although it is not specifically mentioned that the extreme periods are assigned the weight of \(w = 1/N_{i}\). 


\nomenclature{${n}_i$}{number of days after exclusion of extreme days}
\nomenclature{${N}_i$}{number of days}
\nomenclature{$w$}{weight of cluster}

\subsection{Extreme period inclusion methods used in this paper}
\label{sec:Used_methods}

\begin{figure*}[!h]
\centering
{\noindent\includegraphics[width=1.0\textwidth]{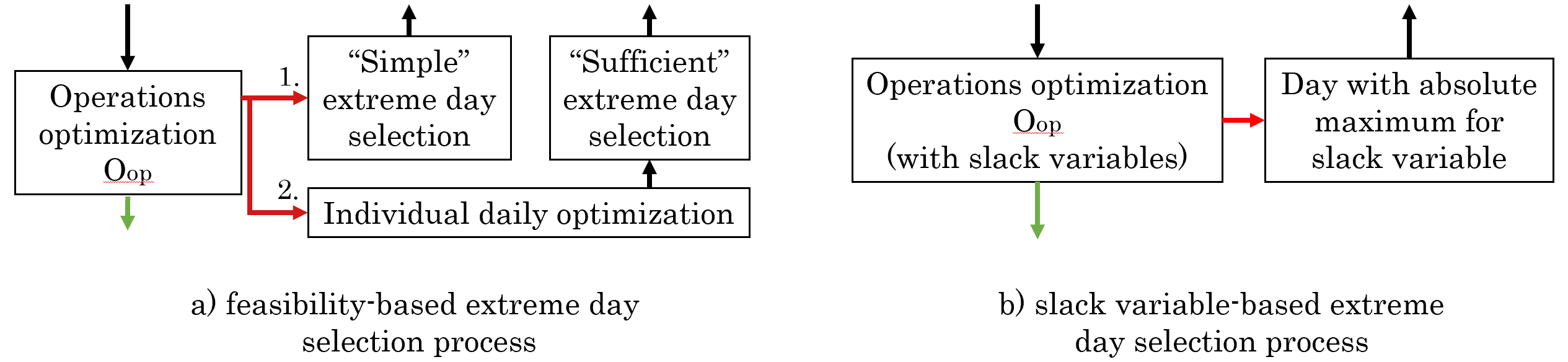}}
\caption{Details on the implementation of the operations optimization and extreme period identification steps of the iterative extreme period selection process (see Figure~\ref{fig:02_02_Feasibility_Based}) to (a) implement the feasibility-based extreme period selection process and (b) the slack variable-based extreme period selection process.}
\label{fig:02_03_Slack_Variable_Based}
\end{figure*}

In this work, we use k=5 clusters obtained from k-means clustering as the basis for our analysis (see Section \ref{sec:Number_Clust} for sensitivity analysis with $k>5$). We investigate three methods to include extreme periods in detail, all using actual periods (An examination of using virtual periods can be found in the SI. We find that virtual periods are too conservative an estimate of the input data distribution). We investigate a simple selection process with statistical identification measure, one hybrid method, and introduce an iterative slack-variable-based method. 
Both latter methods (hybrid and iterative) broadly follow the overall selection process outlined in Figure~\ref{fig:02_02_Feasibility_Based} with some minor modifications. We examine all methods with both ``feasibility'' and ``append'' representation modification. We present the three methods in more detail in the following. 

The first extreme period inclusion method we examine is based on a statistical identification measure and a simple selection process. The extreme periods are then added to the representation obtained from the clustered input data. In our analysis, we have three ``simple'' extreme periods, leading to \(k=5+3\) representative periods. The three ``simple'' extreme periods are the absolute maximum for electricity and heat demand and the integral minimum of the solar availability. 

Which events have to be included as ``simple'' extreme periods depends largely on the characteristics of the optimization problem itself and is up to the modeler's intuition of the problem. For example, for a heat system with a heat pump and electric heater, but without connection to an external heating supply grid, the addition of the peak heat demand would be crucial. There does not need to be an extreme period for all input data sets \(N_{d}\). 

\nomenclature{${N}_d$}{number of input data sets}

The second extreme period inclusion method we examine is a combination of the previously described ``simple'' process and an iterative process with feasibility-based identification measure that adds extreme periods until the convergence criterion that the operations optimization (\(O_{op}\)) is feasible is fulfilled. We call this feasibility-based extreme period selection, and it is adapted from the method developed by Bahl et al.~\citep{Bahl2017a}.
Figure~\ref{fig:02_03_Slack_Variable_Based}a shows how this method fits into the overall selection process that is presented in Figure~\ref{fig:02_02_Feasibility_Based}.
In this method, ``simple'' extreme periods are first identified as described above. By adding ``simple'' extreme periods as a first step, we greatly reduce the number of infeasible days and simplify the iterative process. However, depending on the problem, it is possible that there are additional extreme periods not identifiable a priori by statistical metrics. 

If the operations optimization is infeasible after adding ``simple'' extreme periods, the operations optimization is rerun separately for each day ($O_{daily}$). Thereby, we can determine which specific day or days cause the infeasibility for full operations.
Days that are infeasible for daily operations are then added iteratively one at a time to the set of representative days, until a feasible solution in the operations optimization is reached. We call the combination of ``simple'' and secondary extreme period ``sufficient'' extreme periods (as the combination is sufficient to guarantee feasibility). 

This leads to \(k=5+X\) representative days, where X is the number of ``sufficient'' extreme periods. Depending on the optimization problem at hand, the number of ``sufficient'' extreme periods (\(X\)) varies and thus cannot be stated generally.

Because we add extreme periods successively based on infeasibility, this method does not necessarily ensure that we use the minimum number of extreme periods. Furthermore, the information that a day is infeasible does not say anything about how much it contributes to the infeasibility, or if adding that day will also remove infeasibility that is caused due to other days. 

The third extreme period inclusion method we examine is a newly introduced slack-variable-based identification measure. This method adds extreme periods until all slack variables of the operations optimization (\(O_{op}\)) are zero. We call this method ``slack-variable-based'' extreme period selection. To the best of our knowledge, slack-variable based selection has not previously been introduced in the literature.
Figure~\ref{fig:02_03_Slack_Variable_Based}b shows how this method fits into the overall selection process shown in Figure~\ref{fig:02_02_Feasibility_Based}.

If the operations optimization \(O_{op}\)  is infeasible, it is rerun with slack variables added to the demand constraints (see Section \ref{sec:Opt_problem}).
Slack variables are additional variables that allow the energy system to ``supply'' additional energy in order to achieve feasibility of the system setup. This means that if the system design is underestimated and cannot provide enough energy to meet the demand, the slack variables provide the missing energy at high cost and achieve feasibility. 

During the operations optimization we determine the slack variables for every hour of every day of the original input data set. We identify the day with the absolute maximum of each slack variable as the selected extreme period. This day is then added to the representative days of the clustering process. Then, slack variables of the full operations problem (\(O_{op}\)) are evaluated again and checked if the convergence criterion is reached. If not, the process is iteratively continued until all slack variables of the full operations problem (\(O_{op}\)) are zero. 

Note that for every energy type of the optimization problem (e.g.~electricity, heat), one slack variable needs to be introduced. It is up to the modeler to decide from which energy type the first extreme period should be added, and there is unlikely to be a single approach that works well for all problems. 

In our problem below we use the heat demand slack variable first to identify extreme period, as it is the more constrained energy type. This is due to the lack of supply grid for heat in our optimization problem, whereas power can be supplied via a central grid access (see Section \ref{sec:Opt_problem} for details). 

Note that the slack variable-based process is iterative and that the order in which extreme periods are added may influence the number of necessary extreme periods. Similar to the feasibility-based extreme period selection process the number of ``sufficient'' extreme periods (\(X\)) varies and thus cannot be stated generally.

\nomenclature{X}{unknown number of added extreme days}

\section{Residential energy supply system and input data}
\label{sec:Data}

\subsection{Residential energy supply system}
\label{sec:Opt_problem}

\begin{figure}[!h]
\centering
{\noindent\includegraphics[width=0.5\textwidth]{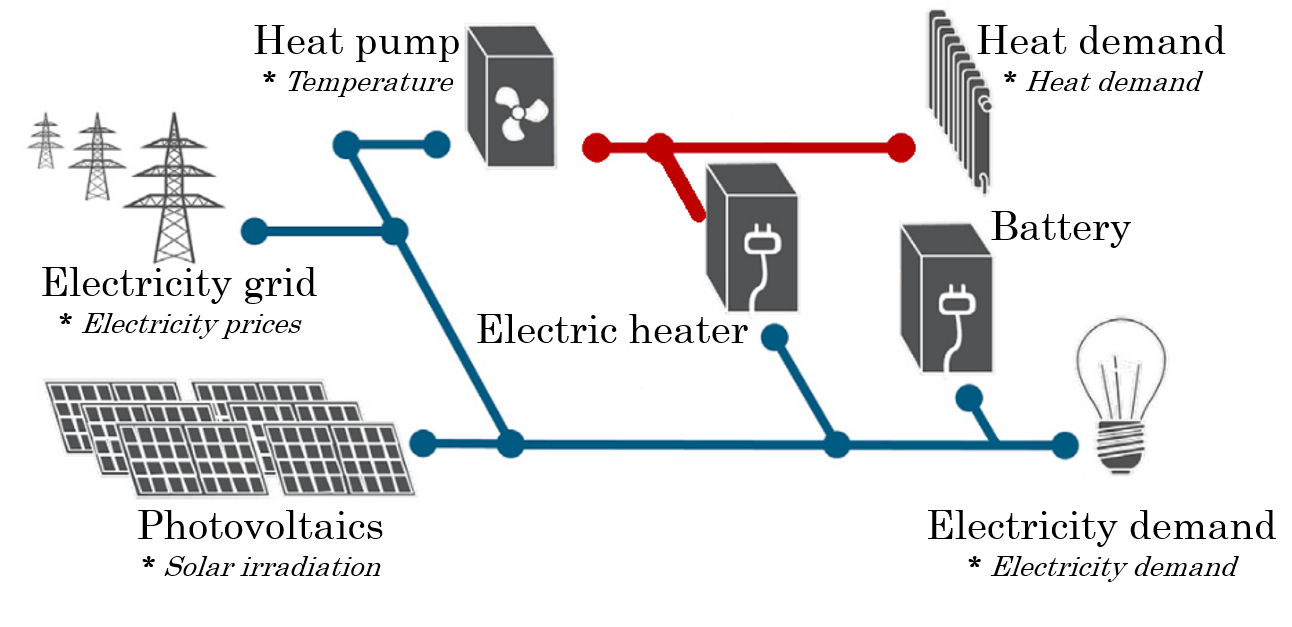}}
\setlength{\abovecaptionskip}{5pt}
\caption[System overview: residential energy supply system with heat and electricity]{System overview: residential energy supply system for electricity and heat. Design variables are battery size, PV plant size, heat pump size, and electric heater size. The asterix indicates the time-varying input data sets that are used in the clustering process. The overall system is adopted from Kotzur et al.~\citep{Kotzur2018}. }
\label{fig:03_01_problem_set}
\end{figure}

We use a modified version of the residential energy supply system introduced by Kotzur et al.~\citep{Kotzur2018}, which was developed to compare time-series aggregation methods. Figure~\ref{fig:03_01_problem_set} shows our modified residential energy supply system to integrate extreme periods and highlights where input data are time-varying and thus clustered. The system consists of a heating system with a heat pump and an electric heater to provide heat via electricity, a photovoltaic system, a central grid connection to provide electricity demand from the central grid, and a battery for storage. 
We modify the system originally presented by Kotzur et al.~\citep{Kotzur2018} in the following ways: First, we use a battery as a storage technology to store electric energy instead of heating storage. Second, we use time-varying electricity prices, which are thus clustered. Third, we limit our grid connection instead of having an unlimited central grid connection. These additions allow for in-detail evaluation of time-varying effects of the electricity attribute, which is also important in many other optimization problems in the literature.

In this work, we use the battery only as an intra-day storage and do not allow for inter-day or seasonal storage. This modeling assumption will likely underestimate the battery usage for a real application, as inter-day storage may reduce overall system cost. However, this modeling assumption affects all calculations and thus does not influence the comparison of the different extreme period methods. For more about intra-day and seasonal storage while using clustered input data we refer the reader to Kotzur et al.~\citep{Kotzur2018a} and Gabrielli et al.~\citep{Gabrielli2018}.

We implement the optimization problem in ``Julia'' \citep{Bezanson2014}, using the ``JuMP'' package \citep{Dunning2015} and ``CPLEX 12.8'' as the solver. 

The optimization problem consists of design and operational decisions variables. The design variables are the sizes of each system part: photovoltaic power output (\(P^{PV}\)) [W], heat pump size (\(P^{HP}\)) [W], electric heater size (\(P^{EH}\)) [W] and battery size (energy capacity (\(E^{bat}\)) [kWh] and battery power capacity(\(P^{bat}\))) [W]. The operational decision variable are the amount of energy purchased from the central electricity grid (\(E^{buy}_{j,t}\)) [kWh], the battery charge (\(E^{in}_{j,t}\)) [kWh] and discharge (\(E^{out}_{j,t}\)) [kWh], the battery storage level (\(Stor^{lev}_{j,t}\)) [kWh], and the electricity used by the electric heater (\(E^{EH_{el}}_{j,t}\)) [kWh] and heat pump (\(E^{HP_{el}}_{j,t}\)) [kWh].

\nomenclature{$P^{PV}$}{photovoltaic plant size}
\nomenclature{$P^{HP}$}{heat pump size}
\nomenclature{$P^{EH}$}{electric heater size}
\nomenclature{$P^{bat}$}{battery power capacity}
\nomenclature{$E^{bat}$}{battery energy capacity}
\nomenclature{}{}

We introduce the following sets for notation: $DV$ $=$ $\{P^{HP}$, $P^{EH}$, $P^{PV}$, $P^{bat}$, $E^{bat}\}$ is the set of design decision variables (the battery is associated with two variables, one for power capacity and one for energy capacity).
The full time-series consists of \(N_{full} = 365\) days, each having \(N_{t} = 24\) hours. \(t \ \epsilon \  T = \{1, ..., N_{t}\}\) is the set of time-steps within one day, in our case hourly. The input data are clustered into \(k\) clusters. \(j \ \epsilon \  K = \{1, ..., k\}\) is the set of cluster indices.

\nomenclature{v}{set of decision variables}
\nomenclature{t}{set of time-steps within one day}
\nomenclature{j}{set of clusters}

We model the problem as a linear program. This allows for \(O_{ref}\) to be computationally tractable. The objective function represents the overall cost of the residential energy supply system and is calculated as follows:

\begin{equation}
\min \left[  \sum_{j=1}^{k} N_{j}* \sum_{t=1}^{N_{t}} \left( E^{buy}_{j,t}*\tilde{c}^{el}_{j,t} \right) + \sum_{v\epsilon V} APVF_{v} * DV_{v} * c_{capex,v} \right]
\label{for:03_01_01-obj}
\end{equation}
\nomenclature{APVF}{annuity present value factor}
\nomenclature{$\tilde{c}^{el}$}{clustered, hourly electricity price}
where \(N_{j}\) is the number of days assigned to each particular cluster \(k\). \(E^{buy}_{j,t}\) is the amount of electricity bought from the central electricity grid at any given hour \(t\) for any given representative day \(j\). \(\tilde{c}^{el}_{j,t}\) represents the day-ahead electricity price for each hour for every representative day. All variables featured with a tilde are input data sets which are clustered. \(APVF_{v}\) is the annuity present value factor and \(c_{capex,v}\) states the power-specific and energy-specific capital cost for each design variable. \(DV_{v}\) is the size of each design variable (\(P^{PV}\), \(P^{HP}\), \(P^{EH}\), \(E^{bat}\), and \(P^{bat}\)).
The capital cost can be found in the SI. 

The operating cost in our system consists only of the cost of buying electricity from the central grid. All other operation costs are excluded (for example for maintenance), from the objective function, as those cost are small compared to the capital and electricity purchasing operational cost. 

In order to investigate the system under different constraining conditions, we limit the central grid connection defined by  \(C_{lim}\) for our system in certain cases: The maximum power supply through the central grid connection necessary to support the system for the full input data optimization ($O_{ref}$) is 5.84 kW. We call this ``100\% grid connection'', which equals \(C_{lim} = 5.84 \ kW\). We define the central grid connection by the maximum power it can transfer. For our experiments, we change the grid connection from 120\% to 0\% (\(C_{lim,\ max} = 7.00\  kW, C_{lim,\ min} = 0 \ kW\)), thus creating different self-sufficiency levels for our system. 

\nomenclature{$C_{lim}$}{cable limit}

In cases that include slack variables, the objective function, originally Equation \ref{for:03_01_01-obj}, is modified as follows:

\begin{multline}
	 \min   \ \sum_{j=1}^{k} N_{k}* \sum_{t=1}^{N_{t}} ( E^{buy}_{j,t}*\tilde{c}^{el}_{j,t} + (E^{slack,\ el}_{j,t} + Q^{slack,\ heat}_{j,t})*c_{slack} ) \\ + \sum_{v\in V} APVF_{v} * DV_{v} * c_{capex,v}
\label{for:03_01_03-obj_slack}
\end{multline}
\nomenclature{$E^{slack,\ el}_{j,t}$}{slack variable for electricity}
\nomenclature{$Q^{slack,\ heat}_{j,t}$}{slack variable for heat}
\nomenclature{$c_{slack}$}{cost of energy for slack variables --- VoLL}
where \(E^{slack,\ el}_{j,t}\) and \(Q^{slack,\ heat}_{j,t}\) are the time- and cluster-dependent slack variables of additional electricity and heat. \(c_{slack}\) is the price for this additional energy. In the literature, this price is also called Value of Lost Load (VoLL). 
We add two slack variables because supply shortages can occur in both energy types (electricity and heat).

Besides the objective function, we add the slack variables (\(E^{slack,\ el}_{j,t}\) and \(Q^{slack,\ heat}_{j,t}\)) to the corresponding constraints for electricity and heat supply.




\subsection{Input data}
\label{sec:Input_data}

The following input data sets are clustered for the optimization problem presented in this section: Electricity and heat demand, ambient temperature, solar availability, and electricity price. Please refer to the SI for additional details on input data.
In total, all data sets but electricity prices affect the constraints. Electricity prices do not affect the feasibility of our optimization problem, but occur in the objective function of the optimization problem. 

One particular factor which affects the selection of ``simple'' extreme periods is that the absolute peaks for heat demand and ambient temperature occur simultaneously (highest heat demand and lowest ambient temperature), so that only one extreme period for those two data sets is necessary. In general, this may not be the case, and would result in additional extreme periods.

\section{Results}
\label{sec:Results}

All data sets are clustered using the k-means clustering method as described in Section \ref{sec:Methods}. For a comparison of the clustered input data versus the full input data, and to measure the effectiveness of the clustering on a similar problem, we refer to Kotzur et al.~\citep{Kotzur2018}. 
Because we analyze daily periods in this study, we refer to extreme periods as extreme days in this section. 
Unless noted otherwise, we use \(k=5\) representative days plus extreme days for our analysis. 

Our results are structured in three main points: First, we analyze the general effect of adding extreme days to the clustered input data of our residential energy supply system. Second, we compare the optimization results (design variables and objective function value) of the two introduced extreme day selection processes (feasibility-based and slack variable-based) and both representation modification methods. Third, we compare the effect of adding additional data via additional clusters vs.~via adding extreme days.

\subsection{Effects of adding simple and iteratively selected extreme days}
\label{sec:Effects}

\begin{figure*}[!h]
\centering
{\noindent\includegraphics[width=0.8\textwidth]{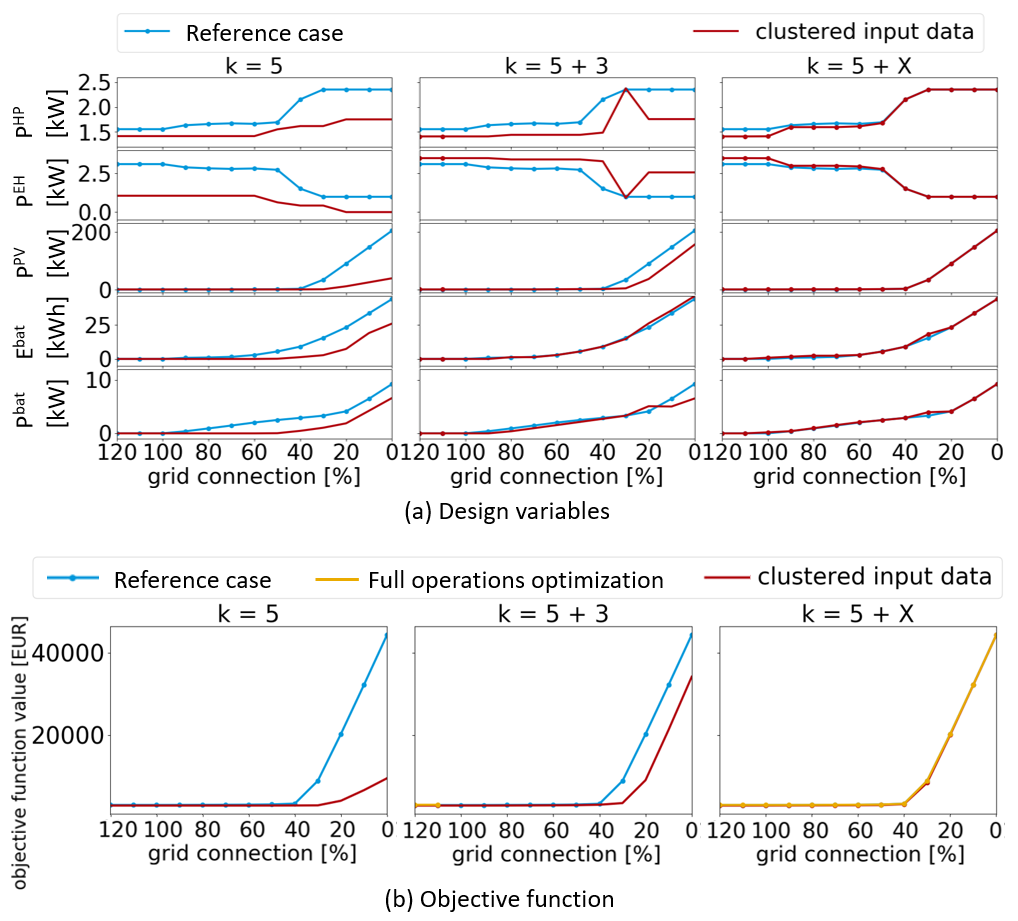}}
\setlength{\abovecaptionskip}{5pt}
\caption[Design variables (a) and objective function value (b) for various grid connections and different extreme day selections]{Design variables (a) and objective function value (b) of the residential energy supply system with no extreme days (k = 5), ``simple'' extreme days (k = 5+3) and ``sufficient'' extreme days (k = 5+X, where X is a number that varies depending on the case that is analyzed) added. The representation modification method is ``feasibility steps''. The grid connection varies from 120\% to 0\%. The input data are clustered in k = 5 clusters plus the possible extreme days. 
From top to bottom we show the capacities of heat pump (\(P^{HP}\)), electric heater (\(P^{EH}\)), photovoltaics (\(P^{PV}\)), energy battery energy (\(E^{bat}\)) and battery power (\(P^{bat}\)). 
The blue line indicates the design variables for the reference system (full input data over 365 days) optimization. The red line shows the design variables for the optimization results with clustered input data. In subplot b the yellow line shows the objective function value of the full operations optimization ($O_{op}$). The dots indicate that a particular system design for a certain grid connection for the full year operations optimization is feasible. }
\label{fig:04_01_DV_k=5_abs}
\end{figure*}

Figure~\ref{fig:04_01_DV_k=5_abs} shows the design variables (Figure~\ref{fig:04_01_DV_k=5_abs}a) and objective function value (Figure~\ref{fig:04_01_DV_k=5_abs}b) for the residential energy supply system for the reference case ($O_{ref}$, full input data) and for the clustered input data with and without the addition of extreme days. Depending on the grid connection, the number of added ``sufficient'' extreme days varies and thus cannot be stated generally. The extreme days are selected by the iterative feasibility-based selection process (see section \ref{sec:Selection}).
Each row of Figure~\ref{fig:04_01_DV_k=5_abs}a illustrates the size of one design variable.

\nomenclature{obj}{objective function value}

The first two rows of Figure~\ref{fig:04_01_DV_k=5_abs}a show the heating design of our system: Without the addition of extreme days, both heat pump and electric heater are underestimated and thus cannot provide enough heat for every day of the year. This the system design infeasible for the 365 day operations optimization. With the addition of ``simple'' extreme days, the heating system is sized larger, though still considerably different from the reference system. The addition of ``sufficient'' extreme days produces reasonably accurate results compared to the reference system.

Similar effects appear for the other design variables. The system design is significantly underestimated without the addition of extreme days. Thus, the system cannot provide energy for the full year, making the full operations optimization $O_{op}$ infeasible. This is especially true for low grid connections (right side of each column). Adding ``simple'' extreme days brings (middle) the design variables closer to a feasible system (e.g. full input data shown in blue). However, it still is not feasible, except for grid connections over 100 \%. Nevertheless, just the addition of ``simple'' extreme days results in a significant increase in accuracy for the design variables, especially for grid connections above 40\%. 

To ensure the feasibility for all grid connections, the addition of ``sufficient'', iteratively selected, extreme days is necessary. This shows once more the importance of testing for feasibility for operations over 365 days, as introduced by Bahl et al.~\citep{Bahl2017a}. In particular for limited grid connections, the addition of extreme days is pivotal. Without the addition of those ``sufficient'' extreme days, the PV plant would be undersized, thus not provide enough power to supply the system for the entire year. 

Figure~\ref{fig:04_01_DV_k=5_abs}b shows the objective function value of the optimization problems in the same structure as Figure~\ref{fig:04_01_DV_k=5_abs}a. Additionally, the yellow line shows the objective function value of the full operations optimization. Please note that the full optimization operations results can only be obtained if the system design is feasible. Overall, Figure~\ref{fig:04_01_DV_k=5_abs}b shows similar results as for the design variables in Figure~\ref{fig:04_01_DV_k=5_abs}a: The addition of extreme days, first ``simple'' and eventually ``sufficient'' extreme days, leads to a higher accuracy in terms of the objective function value, especially for lower grid connections.
However, while the accuracy of the clustered input data objective function value is about 95\% , the accuracy in terms of the full operational optimization is over 99\% compared to $O_{ref}$, see Figure~\ref{fig:04_01_DV_k=5_abs}b.
 
Furthermore, not only the objective function value itself varies, the ratio of operational cost to capital cost changes with varying the grid connection. Table~\ref{tab:04_01-Total_cost_CAPEX_OPEX} shows this change: The share of the capital cost to the total system cost increases with lower grid connections and the share of operational cost decreases with lower grid connections. In a self-sufficient system (0\% grid connection) there is no cost for buying electricity, and therefore the total operational cost of the system is zero. Instead, the objective consists only of minimizing the capital cost, which depends on the extreme cases of the input data. 

\begin{table*}[ht]
\centering
\caption{Overview of total system cost, capital cost, and operational cost for different grid connections for k = 5 + X clusters (``sufficient'' extreme days to ensure feasibility) using feasibility-based extreme day selection process. Note that values are rounded.}
\label{tab:04_01-Total_cost_CAPEX_OPEX}
\begin{tabular}{ccccc}
\toprule
\begin{tabular}[t]{@{}l@{}} \textbf{Grid connection} \\ \textbf{[\%]} \end{tabular}		& \begin{tabular}[t]{@{}l@{}} \textbf{total system cost} \\ \textbf{[EUR/year]} \end{tabular}		& \begin{tabular}[t]{@{}l@{}}  \textbf{capital} \\ \textbf{cost ratio}	\end{tabular}	& \begin{tabular}[t]{@{}l@{}}  \textbf{operational} \\ \textbf{cost ratio} \end{tabular} & \begin{tabular}[t]{@{}l@{}}  \textbf{number of} \\ \textbf{extreme periods (X)} \end{tabular}   \\
\midrule

\addlinespace
120 		&		\ 2,910		& 		0.231		& 		0.769 & 3\\
110 		&		\ 2,910		&		0.231   	& 		0.769 & 3\\
100 		& 		\ 2,914		& 		0.245		& 		0.755 & 4\\
\ 90		&		\ 2,928		&		0.246	    & 		0.754 & 4\\
\ 80    	&		\ 2,948		& 		0.298		& 		0.702 & 4\\
\ 70 		&		\ 2,968		&   	0.303		& 		0.697 & 4\\
\ 60 		&		\ 2,990		& 		0.322		& 		0.678 & 4\\
\ 50 	    &		\ 3,027		&		0.391		& 		0.609 & 4\\
\ 40    	&		\ 3,255		& 		0.484		& 		0.516 & 5\\
\ 30 		&		\ 8,444		&		0.985		& 		0.015 & 5\\
\ 20 		&		20,039		& 		0.999		& 		0.001 & 5\\
\ 10 		&		32,167		&		1.000		& 		0.000 & 5\\
\ \ 0 		&		44,337		& 		1.000		& 		0.000 & 5\\
\bottomrule
\end{tabular}
\end{table*}

Overall, the lower the grid connection, the higher the importance of the capital cost and the lower the importance of the operational cost. The tipping point, where the operational cost becomes significantly less important, occurs with the large increase in design capacity at around 40\% to 30\% grid connection. For grid connections below this point, it appears that the nature of the optimization problem changes fundamentally and the system cannot be supplied by the constrained central grid connection for most of the time. 

This is to be expected, because a lower grid connection directly results in less purchased electricity. It also means that the clustering process for lower grid connections becomes less and less important  for the design of the energy system and the correct identification of extreme days gains relevance. Once the grid connection becomes low enough, the design of the system eventually depends solely on the correct extreme day selection rather than on the clustered input data, which are more average days. 

Overall we find that without the selection and addition of the correct extreme days, the system may not be feasible for  365 day operations and lack accuracy of the objective function, especially with limiting central grid connection or similar conditions. Hence, the challenge lies in identifying these extreme days for each grid connection and, more generally, for each optimization problem. 
For limited grid connections of below 40\%, extreme periods become imperative to ensure reasonable accuracy of design variables and objective function.

\subsection{Comparison of selection and representation modification methods}
\label{sec:Comparison}
 
 \begin{figure*}[!h]
\centering
{\noindent\includegraphics[width=0.7\textwidth]{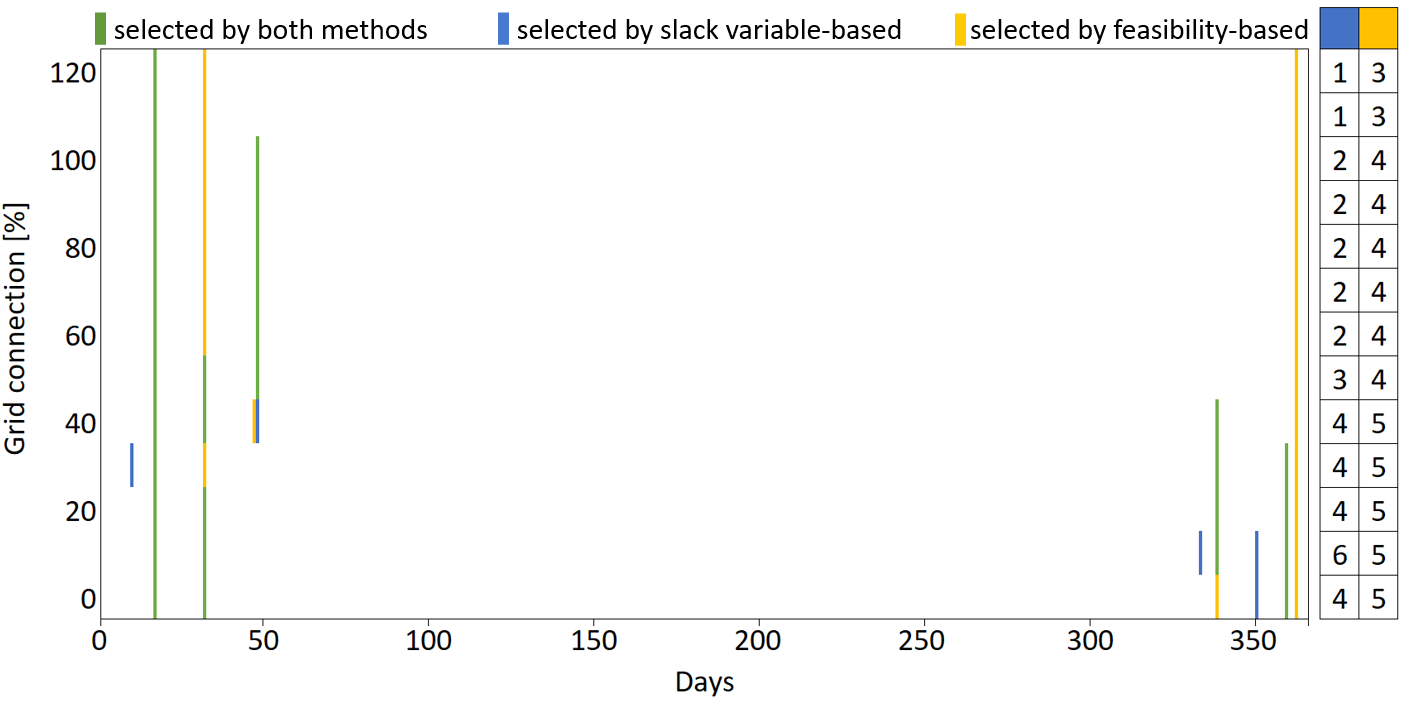}}
\setlength{\abovecaptionskip}{5pt}
\caption{Selected extreme days by the two extreme day selection processes for various grid connections. The days marked in green are selected by both the feasibility-based extreme day selection process and the slack variable-based extreme day selection process, whereas the days marked in yellow are selected by the former process only, and the days marked in blue are selected by the latter process only. Numbers on the right indicate the total number of days that were selected by each method.}
\label{fig:04_02_Selected_Extreme_Days_All_Methods}
\end{figure*}

In this paper we use two different iterative processes which can ensure feasibility of the optimization problem solution: feasibility-based extreme day selection and slack variable-based extreme day selection process and two different modifications for including the selected extreme days to the clustered input data. 

Figure~\ref{fig:04_02_Selected_Extreme_Days_All_Methods} shows the different days each extreme day selection process selected to ensure the feasibility of the full operations optimization problem (\(O_{op}\)) for various grid connections from 120\% to 0\%. In general, the lower the grid connection, the more extreme days need to be added to the clustered input data in order to ensure feasibility, regardless of the chosen selection process. 
Figure~\ref{fig:04_02_Selected_Extreme_Days_All_Methods} shows that the number of selected extreme days is influenced by both the grid connection and whether the selection process is feasibility-based or slack variable based: The slack variable-based selection process selects on average fewer extreme days in order to achieve a feasible optimization result.

\begin{figure}[!h]
\centering
{\noindent\includegraphics[width=0.5\textwidth]{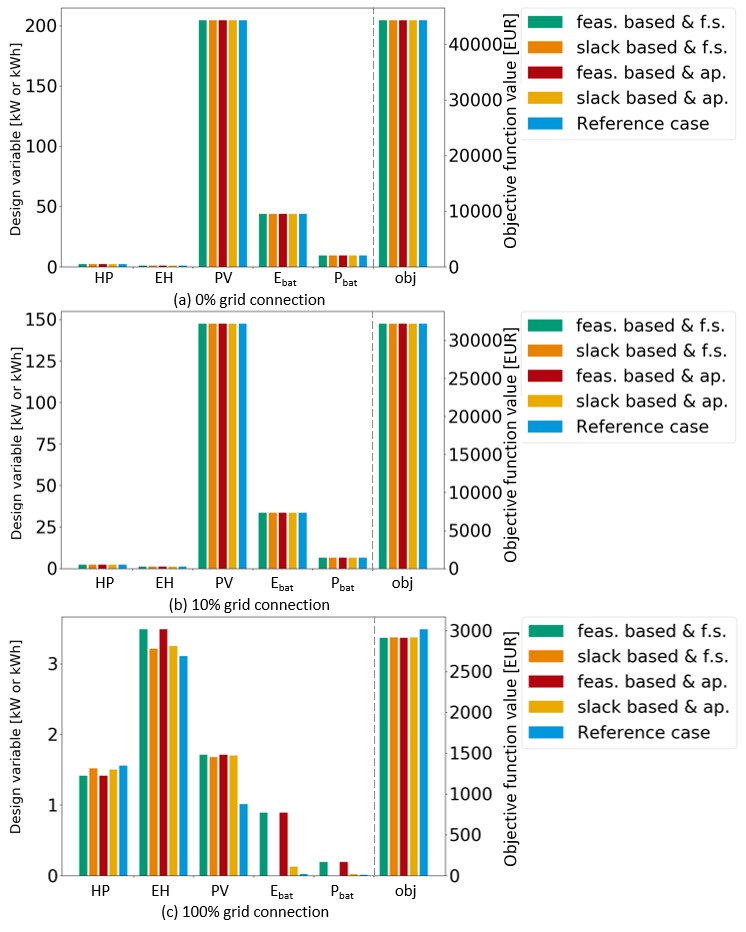}}
\setlength{\abovecaptionskip}{5pt}
\caption{Comparison of design variables and objective function value of the optimization problem for the different combinations of extreme day selection -- feasibility-based extreme day selection and slack variable-based extreme day selection -- and representation modification -- ``f.s.'' (feasibility steps) and ``ap.'' (append) -- for different grid connections. The number of clusters is k=5 plus the identified extreme days shown in Figure~\ref{fig:04_02_Selected_Extreme_Days_All_Methods}. This results in k=5+5 (a), k=5+5 (b) and k=5+4 (c) for the feasibility-based selection and k=5+4 (a), k=5+6 (b) and k=5+2 (c) for the slack variable-based selection.}
\label{fig:04_03_Comparision_Selection_Representation}
\end{figure}

Figure~\ref{fig:04_03_Comparision_Selection_Representation} shows the effect of selected extreme days via the two different selection processes on the results of the optimization problem, while using both representation modification methods, ``feasibility steps'' and ``append''. Figure~\ref{fig:04_03_Comparision_Selection_Representation}a shows the comparison of all five cases for a grid connection of 0\%, Figure~\ref{fig:04_03_Comparision_Selection_Representation}b for 10\% and Figure~\ref{fig:04_03_Comparision_Selection_Representation}c for 100\%. Those three grid connections are chosen, due to the different number of extreme days that were selected by the two extreme day selection processes.

Figure~\ref{fig:04_03_Comparision_Selection_Representation} illustrates that there is little effect on the design variables or objective function of the optimization problem, even though different extreme days have been selected and both the number and specific identification of individual days as extreme days deviates as shown in Figure~\ref{fig:04_02_Selected_Extreme_Days_All_Methods}.

This means that it is not necessary to identify the exact critical extreme days (minimal days to ensure feasibility), but that it is possible to achieve very similar accuracy with the addition of more extreme days than necessary. This is true, as long as all infeasibility causing days are identified eventually. However, the lower the number of extreme days that need to be added, the better. 

Moreover, Figure~\ref{fig:04_03_Comparision_Selection_Representation} shows the accuracy for design variables and for the objective function value of the optimization problem with both representation modification methods, ``feasibility steps'' and ``append''.
For 0\% grid connection there are no differences between the different representation modification methods, regardless of the extreme day selection process. This is to be expected, because only the infeasibility causing days determine the design for this case, as explained in Section \ref{sec:Effects}. For a grid connection of 10\%, similar effects are shown. The different extreme days do not lead to any significant changes in the optimization results, neither for the design variables nor the objective function value. For a grid connection of 100\% the results of the design variables and objective function value of the optimization problem vary slightly, however, this variation does not depend on the selected extreme days (different extreme day selection processes). Thus, for different selected extreme days, the two representation modification methods result in almost identical results, leaving the choice of extreme day selection process open. Furthermore, the changes in objective function value caused by the different representation modification methods are insignificantly small (\(< 0.1\%\)). 

Overall, Figure~\ref{fig:04_03_Comparision_Selection_Representation} shows that both representation modification methods perform equally well for both extreme day selection processes, even when different extreme days are selected. However, because the ``append'' method requires recomputation of the clusters, we conclude that ``feasibility steps'' are the simplest and yet very accurate representation modification to use on our problem. 
Moreover, for the implementation of iterative processes,  ``feasibility steps'' are computationally less expensive because they do not need the clustering to be performed again after each identified extreme. 

\subsection{Effects of number of clusters}
\label{sec:Number_Clust}

Thus far, we have compared the different identification, selection and representation of extreme days with k=5 clusters. Previous analyses of representative days selected by the k-means algorithm showed that more representative days can lead to an increase in accuracy of the optimization results \citep{Nahmmacher2016, Bahl2017a, Kotzur2018, Teichgraeber}. Now we vary k and compare the influence of adding additional clusters instead of adding extreme days instead. We make this comparison for the grid connection of 0\% and 100\%. In both cases four days have been added as sufficient extreme days to ensure feasibility.

\begin{figure}[h]
\centering
{\noindent\includegraphics[width=0.5\textwidth]{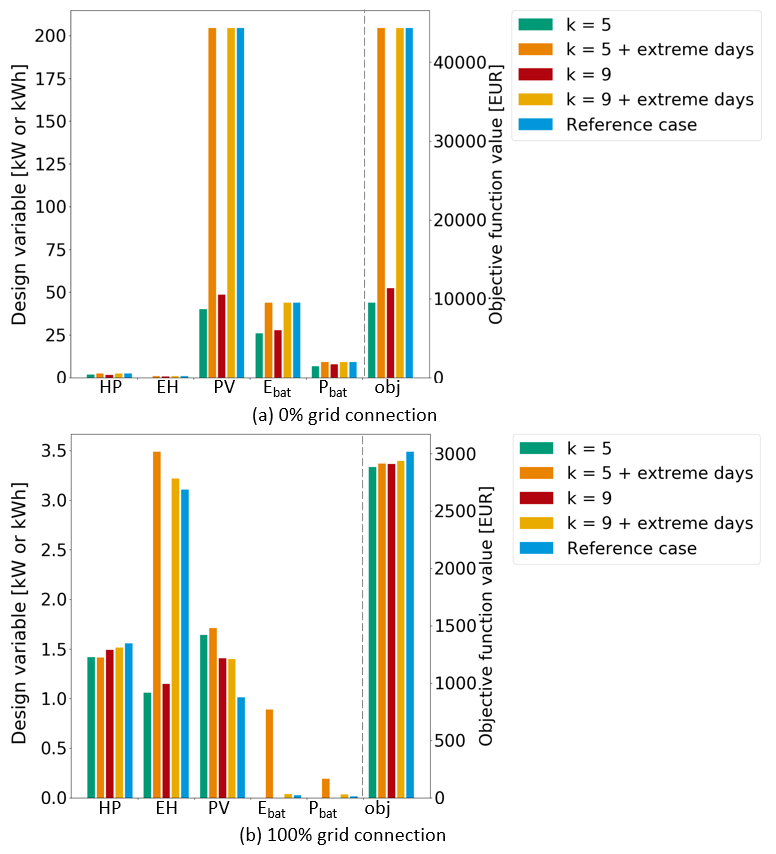}}
\setlength{\abovecaptionskip}{5pt}
\caption[Comparison of various number of clusters with and without extreme days]{Comparison of k=5 and k=9 clusters with and without sufficient extreme days selected by the feasibility-based extreme day selection process for all design variables and the objective function of the optimization problem. We show the comparison for (a) 0\% grid connection and for (b) 100\% grid connection. }
\label{fig:04_05_k=5_k=9}
\end{figure}

Figure~\ref{fig:04_05_k=5_k=9} shows all design variables and the objective function of the optimization problem for k=5 clusters and k=9 clusters with and without ``sufficient'' extreme days (k=5+4 and k=9+4 respectively). For a grid connection of 0\%, both design variables and the objective function are significantly underestimated without the addition of extreme days, regardless of whether 5 or 9 representative days were used, although k=9 clusters results in slightly better accuracy than k=5 clusters. With the addition of the sufficient extreme days, both optimization runs, with k=5 clusters plus extreme days (k=5+4) and k=9 clusters plus extreme days (k=9+4), result in the exact same system design and total system cost. Again this is to be expected, as for 0\% grid connection the design of the energy system relies purely on build capacity and thus extreme days. 

For 100\% grid connection, there are minor differences in terms of design variables and objective function with different numbers of clusters. Note that the system design is only feasible with added extreme periods added. See SI for more details on the 100\% grid connection case. 

For heavily constrained optimization problems (e.g. 0\% grid connection), adding extreme periods is far more important in order to reach a high accuracy for design variables and total cost estimates than increasing the number of clusters. We conclude that adding ``sufficient'' extreme days not only makes the optimization problem feasible for the full input data but also allows us to use fewer clusters than normally necessary and still account for high accuracy. In our case, just five clusters and four ``sufficient'' extreme days, led to accuracies of over 98\% compared to the objective function value of the reference system while being a feasible system design.

\section{Conclusion}
\label{sec:Conclusion}

Overall, this work provides a systematic framework to use extreme events as part of representative periods for energy systems optimization problems. Our framework allows for intercomparibility of past and future work. 
Furthermore, we introduce an extreme period inclusion method that is based on slack variables of the optimization problem itself. 
We evaluate a variety of extreme period inclusion methods on a case study: a residential energy supply system. 

There are four specific results of this work:
\begin{enumerate}

\item The inclusion of extreme periods is vital to account for high accuracy and practical applications for clustering in optimization problems. While only using clustering can provide optimization results with decent accuracy, we show that neglecting extreme events usually leads to infeasible solutions on the operations optimization problem with full input data, in particular for heavily constrained optimization problems such as with a limited grid connection or high self-sufficiency constraints. Adding peaks heuristically can reduce the number of the days that are infeasible, but only methods validating the design's performance on the full input data lead to reliable system designs and thus optimization results. 

\item We do not identify any significant differences in terms of the optimization results between the two extreme period inclusion methods we evaluate on our sample problem (including extreme days based on infeasible operations on the original data, and based on slack variables on the demand constraints on the original data). That result prevails even if the methods select different extreme periods. Therefore, we conclude that the specific selection of extreme periods is---on our residential energy supply system problem---less important than the importance of ensuring feasibility by use of extreme periods in the first place. 

\item We could find only marginal differences between representing extreme periods with weight zero (``feasibility steps'') or weight one (``append'') in the objective function. The marginal improvements in optimization results by the append method compared to the feasibility steps method are offset by the significantly more computationally expensive calculations for the append method, because the append method requires clustering after every addition of extreme events. 

\item We show that the inclusion of additional extreme periods to the set of representative periods, instead of additional cluster centers, leads to even higher accuracy while achieving feasibility at the same time. 
We therefore suggest to first start with a relatively small number of clusters and find the ``sufficient'' extreme periods to achieve feasibility and to then increase the number of clusters if a higher accuracy in terms of the optimization results is necessary. 
Particularly for low central grid connections or other constrained systems, extreme periods have a much higher impact on the accuracy than an increase in the number of averaged clusters.
\end{enumerate}

There are several avenues for future research. Investigating optimization problems in which constraints affect the search for extreme values differently is one of them. The residential energy supply system problem at hand is constrained by the grid connection, which applies to each hour individually. Thus, this constraint does not link among periods. However, there are other optimization problems in the literature that are constrained by a maximum emissions limit, which is the integral over the emissions over all periods. This constraint does not allow for the separation of days as in the problem presented here and thus the methods evaluated here are not applicable without modification. Moreover, our introduced method could be tested on a spatially separated problem set with different extremes at different time periods.

\section*{Acknowledgment}
This work was supported by the Wells Family Stanford Graduate Fellowship and the Precourt Institute for Energy Seed Grant for HT. 
LK, DS and MR acknowledge the financial support by the Federal Ministry of Economic Affairs and Energy of Germany in the project METIS (project number 03ET4046). 

\printnomenclature


\section*{References}
\bibliography{Mendeley_modified}

\end{document}


\begin{frontmatter}

\title{Supplementary Information - Extreme events in time series aggregation: A case study for optimal residential energy supply systems}

\author[ere]{Holger Teichgraeber \tnoteref{t1}}
\ead{hteich@stanford.edu}
\author[ere]{Constantin P. Lindenmeyer \tnoteref{t1}}
\ead{constantin.lindenmeyer@rwth-aachen.de}
\author[ltt]{Nils Baumg\"artner}
\author[fzj]{Leander Kotzur}
\author[fzj,cfc]{Detlef~Stolten}
\author[fzj]{Martin Robinius}
\author[ltt,jul]{Andr\'{e} Bardow}
\author[ere]{Adam R. Brandt \corref{cor}}
\ead{abrandt@stanford.edu}

\cortext[cor]{Corresponding author. Tel: 650-724-8251}

\tnotetext[t1]{HT and CPL contributed equally to this work.}

\address[ere]{Department of Energy Resources Engineering, Stanford University,
  Green Earth Sciences Building 065, 367 Panama St., Stanford, California, USA}
  
\address[ltt]{Institute of Technical Thermodynamics, RWTH Aachen University, Schinkelstrasse 8, 52062 Aachen, Germany}

\address[fzj]{Institute of Electrochemical Process Engineering (IEK-3), Forschungszentrum J\"ulich GmbH, Wilhelm-Johnen-Str. 52428 J\"ulich, Germany}

\address[jul]{Institute of Energy and Climate Research - Energy Systems Engineering (IEK-10), Forschungszentrum J\"ulich GmbH, 52428 Jü\"ulich, Germany}

\address[cfc]{Chair for Fuel Cells, RWTH Aachen University, c/o Institute of Electrochemical Process Engineering (IEK-3), Forschungszentrum J\"ulich GmbH, Wilhelm-Johnen-Str., 52428 J\"ulich, Germany}

\end{frontmatter}


\section{Introduction}
This document provides additional details concerning the paper ``Extreme events in time series aggregation: A case study for optimal residential energy supply systems''.

\section{Capital cost}
\label{CAPEX}

Table \ref{tab:06_01-CAPEX} shows our assumptions of the capital cost for all design variables. We choose a amortization period of 5 years.   

\begin{table}[h]
\centering
\caption{Overview of capital cost and amortization period of each design variable for our residential energy supply system. We assume an amortization period of 5 years for all design variables.}
\label{tab:06_01-CAPEX}
\begin{tabular}{llrl}
\toprule
\textbf{Design variable} & \textbf{Symbol}	&	\textbf{Capital cost} 	&	\textbf{Unit} \\

\midrule
photovoltaics		&		\(P^{PV}\)	&		900 	\ \ 		&			\(EUR/kW\)  \\
heat pump			&		\(P^{HP}\)	& 			900 	\ \ 		&			\(EUR/kW\)  \\
electric heater			&		\(P^{EH}\)	&		\ \ 50 	\ \ 			&			\(EUR/kW\)  \\
\begin{tabular}[t]{@{}l@{}} battery	power \\ battery energy \end{tabular}			&		\begin{tabular}[t]{@{}l@{}} \(P^{bat}\) \ \   \\ \(E^{bat}\) \ \  \end{tabular}	&			\begin{tabular}[t]{@{}l@{}} 150 \ \   \\ 100 \ \  \end{tabular}  		&		\begin{tabular}[t]{@{}l@{}} \(EUR/kW\) \ \   \\ \(EUR/kWh\) \ \  \end{tabular}  \\
\bottomrule
\end{tabular}
\end{table}

\section{Input data}
An additional note on the input data used in the paper:
The first four input data sets are the same as in the residential energy supply system from Kotzur et al.~\cite{Kotzur2018}. Electricity and heat demand are for a single family household with four members. Additionally, we use time-varying electricity prices. The electricity prices are based on the German day-ahead market, used in a study by Teichgraeber and Brandt \cite{Teichgraeber2019}. 
The electricity prices are customized to a standard residential supply contract. Thus, the mean electricity price is at 0.301 EUR/kWh, which is the average cost for households in Germany. The maximum electricity price is at 0.370 EUR/kWh, the minimum electricity price at 0.190 EUR/kWh. While time-varying electricity prices may not be in place today, they may be indicative of future energy systems. 

\section{Representation modification}
An additional note on representation modification methods from the literature:
Kotzur et al.~\citep{Kotzur2018} use two more representation modification methods: the first method adds a new cluster with the extreme period as its representation (``add cluster''). The second switches the representation of the cluster which includes the identified extreme periods and changes it from the centroid to the extreme period (``move cluster''). Without additional adaptations both methods are not preserving the mean of the original data. Additionally, the  ``move cluster'' method is constrained by the possible number of extreme periods that can be added, as they need to be fewer than the number of original clusters. In this work, we only present two methods (``feasibility steps'' and ``append'') that are generally applicable regardless of the chosen clustering algorithm or the number of clusters and are preserving the mean of the original data. 

\section{Extreme period inclusion methods used in this paper}
An additional note on feasibility-based extreme day selection:
A method could be developed (not applied here) that exhaustively searches through each extreme periods data set and tests each individually, then choosing the set that leads to the minimum number of extreme periods that satisfies operations feasibility for all days. For complex problems, this would likely be computationally expensive and thus may not be practical.

\section{Effects of number of clusters}
A note adding detailed analysis on the 100\% grid connection result: For 100\% grid connection, there are slight differences between the optimal design variables with different numbers of clusters. While the overall size of the heating system design (heat pump and electric heater) is underestimated without the addition of extreme days, the number of clusters slightly changes the ratio between heat pump design and electric heater design. More clusters and thus more data variability push this ratio closer to the optimal system design of the reference system. Additionally, more clusters seem to better depict the solar availability, because the PV design with k=9 clusters is more accurate than the design with k=5 clusters. In this instance, the addition of extreme days seems to have a small negative impact on the accuracy of an individual design variable, compared to the optimal system design of the reference system. However, looking at the objective function of the optimization problem, the total system cost of k=5 clusters plus the four sufficient extreme days (therefore 9 representative days in total) predicts the total cost with a slightly better accuracy than k=9 clusters without any extreme days. Please note, that the system design for k=9 clusters (without extreme days) is not feasible. Using k=9 clusters plus extreme days leads to another slight improvement and to feasibility, but this still shows that the correct identification of extreme days and their addition to the clustered input data leads to as much or more accuracy improvement as using more clusters. 

\section{Virtual extreme days}
\label{Smarch}
When using actual periods (which is the default in all past studies in the literature), the specific period for all input data sets has been selected and added to the clustered data, therefore adding the real-time extreme event. Instead of using periods where all data types are connected by their time of occurence, one could use ``virtual days''. ``Virtual days'' are the combination of extreme events for different input data sets from different days or periods. For example, a ``virtual day'' could include the peak demand for heat, the peak demand for electricity and the minimum daily solar availability, even though all three events occur during different days over the year. Thus an artificial period is created, which does not occur in the real data. 

In our analysis presented in Section 4, we use actual extreme days to add to the clustered input data. However, it is also possible to merge different extremes from the individual input data sets to one ´´virtual'' extreme day. In this case fewer extreme days are necessary, as ``virtual'' extreme days create a worse combination of individual extremes from the input data set, than naturally occurs in the real input data. In our case, two ``virtual'' extreme days are sufficient to always ensure feasibility. 

\begin{figure}[!htb]
\centering
{\noindent\includegraphics[width=0.95\textwidth]{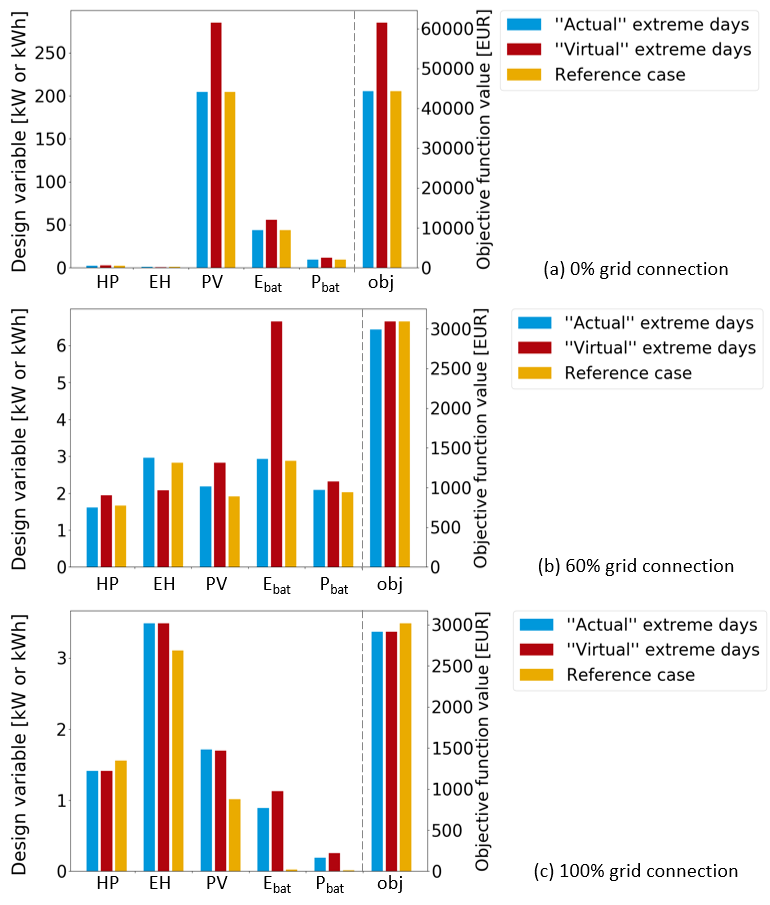}}
\setlength{\abovecaptionskip}{5pt}
\caption[Comparison of ``actual'' and ``virtual'' extreme days]{Comparison of ``actual'' and ``virtual'' extreme days for all design variables and the objective function of the optimization problem. The representation modification method is ``feasibility steps''. We show both extreme day selection methods and the reference case (full input data optimization) for (a) a grid connection of 0\%, for (b) a grid connection of 60\% and for (c) a grid connection of 100\%. }
\label{fig:06_01_virtual_days}
\end{figure}

Figure \ref{fig:06_01_virtual_days} shows the comparison between the use of actual and ``virtual'' extreme days on the optimization results for various grid connections of 0\%, 60\% and 100\%. For 0\% grid connection, ``virtual'' extreme days lead to a massive oversizing resulting from the more extreme peaks of the ``virtual'' extreme day. For a grid connection of 60\% and 100\% there are considerable differences in the sizing of the design variables, however in terms of the objective function value, the resulting differences are relatively small. In case of the 60\% grid connection, the use of ``virtual'' extreme days instead of actual extreme days even leads to a more accurate objective function value compared to the reference case. 

Overall, ``virtual'' extreme days seem to increase the importance of storage technologies, as the battery is most affected by changes when using ``virtual'' extreme days instead of ``actual'' extreme days. However, for lower grid connections the possible effects of overestimation have to be accounted and mitigated for. 

\section*{References}
\bibliography{Mendeley_modified}